\title{On the simple normality to base 2 of $\sqrt{s}$, for $s$ not a perfect square}
\author{Richard Isaac\\
Lehman College and Graduate Center, CUNY}
\documentclass[11pt]{article}

\newtheorem{theorem}{Theorem}

\newtheorem{lemma}{Lemma}

\date{   }
\begin{document}
\maketitle
\begin{abstract}
Since E. Borel proved in 1909 that almost all real numbers with respect
to Lebesgue measure are normal to all bases, an open problem has been whether simple irrational numbers like $\sqrt{2}$ are normal to any base. This paper shows that each number of the
form $\sqrt{s}$ for $s$ not a perfect square is simply normal to the base 2. The argument uses some elementary ideas in the calculus of finite differences.\\
\end{abstract}
\footnote{{\it  \hspace{-2em} AMS 2000 subject classifications}. 11K16. }
\section{Introduction}
A number is {\it simply normal to base $b$} if
its base $b$ expansion has each digit appearing with average
frequency tending to $b^{-1}$. It is {\it normal to base $b$} if its
base $b$ expansion has each block of $n$ digits
appearing with average frequency tending to $b^{-n}$. A number is called
{\it normal}
if it is normal to base $b$ for every base. For a more detailed
introductory discussion we refer to chapter 8 of \cite{ni:ir}.

The most important theorem about normal numbers is the celebrated result
(1909) of E. Borel in which he proved the normality of almost all
numbers with respect to Lebesgue measure. This left open the question,
however, of identifying specific numbers as normal, or even
exhibiting a common irrational normal number. Recently progress has been made in defining certain classes of numbers which can be proved to be normal (see \cite{bc:nn}) but they do not include simple irrationals like $\sqrt{2}$, for example. The difficulty in exhibiting normality for such common irrational numbers is not
surprising since normality is a property depending on the tail of the
base $b$ expansion, that is, on all but a finite number of digits. By
contrast, we mostly ``know'' these numbers by finite approximations, the
complement of the tail.

Identification of well-known irrational numbers as normal may have
interest for computer scientists. The $b$-adic expansion of a number
normal to the base $b$ is a sequence of digits with many of the properties of a random number table. There is thus the possibility that such numbers could be used for the generation of
random numbers for the computer. This would only be possible if the
digits of the normal expansion could be generated quickly enough or
stored efficiently enough to make the method practical.

In this paper we exhibit a class of numbers simply normal to the base 2.
More precisely, we prove
\begin{theorem} \label{snb2} Let $s$ be a natural number which is not a perfect
square. Then the dyadic (base 2) expansion of $\sqrt{s}$ is simply normal.  \end{theorem}

Consider numbers $\omega$ in the unit interval, and represent the dyadic expansion of $\omega$ as
\begin{equation}
\omega=.x_{1}x_{2} \cdots ,\hspace{.4in} x_{i}=0 \mbox{ or } 1.  \label{expom}
\end{equation}
Also of interest is the dyadic expansion of $\nu=\omega^{2}$:
\begin{equation} \label{xu}
 \nu=\omega^{2}=.u_{1}u_{2} \cdots, \hspace{.4in} u_{i}=0 \mbox{ or } 1. \end{equation}
Throughout this paper it will be assumed that $\nu$ is irrational. Then $\omega$ is also irrational and both expansions are uniquely defined. Sometimes it will be convenient to refer to the expansion of $\omega$ as an $X$ sequence and the expansion of $\nu$ as a $U$ sequence. Define the coordinate functions
$X_{n}(\omega)=x_{n}$ and $U_{n}(\omega)=u_{n}$ to be the 
$n$th coordinate of $\omega$ and $\nu$, respectively. We sometimes denote a point of the unit interval by its coordinate representation, that is, $\omega=(x_{1}, x_{2}, \cdots)$ or $\nu=(u_{1}, u_{2}, \cdots)$. Given any dyadic expansion $.s_{1}s_{2}\cdots$ and any positive integer $n$, the sequence of digits $s_{n}, s_{n+1}, \cdots$ is called a {\it tail} of the expansion. Two expansions are said to have the same tail if there exists $n$ so large that the tails of the sequences from the $n$th digit are equal.

The average
\begin{equation} \label{aver}
f_{n}(\omega)= \frac{X_{1}(\omega)+X_{2}(\omega)+\cdots+X_{n}(\omega)}{n} \end{equation}
is the relative frequency of 1's in the first $n$ digits of the expansion of $\omega$. Simple normality for $\omega$ is the assertion that $f_{n}(\omega) \rightarrow 1/2$ as $n$ tends to infinity. Let $n_{i}$ be any fixed subsequence and define \begin{equation} \label{ff}
f(\omega)= \limsup_{i \rightarrow \infty} f_{n_{i}}(\omega). \end{equation}
We note that  
the function $f$ is a tail function with respect to the $X$ sequence, that is, $f(\omega)$ is determined by any tail $x_{n}, x_{n+1}, \cdots$ of its coordinates. In fact, $f$ satisfies a more stringent requirement: it is an invariant function (with respect to the $X$ sequence) in the following sense: let $T$ be the 1-step shift transformation on $\Omega$ to itself given by
\[T(.x_{1}x_{2} \cdots) =.x_{2}x_{3} \cdots. \]
A function $g$ on $\Omega$ is invariant if $g(T\omega)=g(\omega)$ for all $\omega$. Any invariant function is a tail function. 

The following observation will be of particular interest in the proof of Theorem~\ref{snb2}: the average $f_{n}$, defined in terms of the $X$ sequence, can also be expressed as a function $h_{n}(\nu)$ of the $U$ sequence because the $X$ and $U$ sequences uniquely determine each other. This relationship has the simple form  
$f_{n}(\omega)=f_{n}(\sqrt{\nu})=h_{n}(\nu)$. A similar statement holds for any limit function $f$ in relation~\ref{ff}.\\
\noindent {\bf Definition}: Let $f$ be defined as in relation~\ref{ff} for any fixed subsequence $n_{i}$.  
We say that Condition (TU) is satisfied if $f(\omega)=h(\nu)$ is a tail function {\it with respect to the $U$ sequence} whatever the sequence $n_{i}$, that is, for any $\omega$ and any positive integer $n$,  $f(\omega)$ only depends on $u_{n}, u_{n+1}, \cdots$, a tail of the expansion of $\nu=\omega^{2}$.
(The notation ``TU'' is meant to suggest the phrase ``tail with respect to the $U$ sequence''.)
\section{Condition (TU) implies simple normality}
In this section we prove that Condition (TU) implies Theorem~\ref{snb2}.
\begin{theorem}\label{w12} If Condition (TU) is satisfied then Theorem~\ref{snb2} is true.\end{theorem}
The proof requires the following result:
\begin{lemma}   \label{pts2} Let $s$ be a natural number which is not a perfect
square, and let $l$ be any integer such that $2^{l} > s$.
Define the points
\[\omega_{s_{1}}=1-(\sqrt{s}/2^{2l})  \]
and
\[ \omega_{s_{2}}=(\sqrt{s}-1)/2^{l}.  \]
Let $f=\limsup_{i \rightarrow \infty}f_{n_{i}}$ where $n_{i}$ is any fixed subsequence. Assume Condition (TU) is satisfied. Then 
$f(\omega_{s_{1}})=f(\omega_{s_{2}})$. \end{lemma}
Proof: The numbers $\omega_{s_{i}}$ are less than 1 for
$i=1, 2$ and their squares are both irrational and are respectively given by
\begin{eqnarray} \label{twopts}
 1+s(2^{-4l})-(2^{-2l+1} \sqrt{s}) \mbox{  and  }
(s+1)2^{-2l}-(2^{-2l+1} \sqrt{s}). \end{eqnarray}
The dyadic expansions of the rational terms $1+s(2^{-4l})$ and $(s+1)2^{-2l}$
in relation~\ref{twopts} have only a finite number of non-zero digits.
Now consider the dyadic expansion of the term $2^{-2l+1} \sqrt{s}$\, (this
is obtained from the expansion of
$\sqrt{s}$ by shifting the ``decimal'' point $2^{2l-1}$ places to the
left). To get each of the values in relation~\ref{twopts}, this term
must be subtracted from each of the larger rational terms which
have terminating expansions; it is clear that the
resulting numbers have expansions with the same tail, that is, the expansions
of $\omega^{2}_{s_{1}}$ and $\omega^{2}_{s_{2}}$ have the same tail.
Then Condition (TU) implies that
$f(\omega_{s_{1}})=f(\omega_{s_{2}})$. This finishes the proof of the Lemma.

The proof of Theorem~\ref{w12} will now be completed. It is sufficient to prove simple normality for $\lambda=\sqrt{s}-[\sqrt{s}\,] < 1$ where $[t]=$ greatest integer $ \leq t$. Define $g_{n}(\omega)$ to be the average number of 0's in the first $n$ digits of the expansion of $\omega$.
Let $n_{i}$ be any subsequence
such that $f_{n_{i}}(\lambda)$ converges to some value $a$. Now consider
the point $\lambda^{\prime}=1-\lambda$, and notice that for all $j$ the $j$th digit
of $\lambda$ and the $j$th digit of $\lambda^{\prime}$ add to 1. It
follows that $g_{n_{i}}(\lambda^{\prime})$ also converges to $a$. Note
that the point $\omega_{s_{1}}$ (as defined in Lemma~\ref{pts2})
would have the same tail as $\lambda^{\prime}$ were we to shift a finite
number of places, and therefore $\omega_{s_{1}}$ and $\lambda^{\prime}$
have the same asymptotic relative frequency of 0's and 1's. The same can
be asserted for $\omega_{s_{2}}$ and $\lambda$. Thus
Lemma~\ref{pts2} can be applied to conclude that the asymptotic
averages based on $f_{n_{i}}$ evaluated at the points $\lambda$ and
$\lambda^{\prime}$ are equal, that is, 
\[ \limsup_{i \rightarrow \infty}
f_{n_{i}}(\lambda^{\prime} )=\lim_{i \rightarrow \infty} f_{n_{i}}(\lambda) = a . \] But the
equation $f_{n}+g_{n}=1$ holds for all $n$ at all points; apply it for
$n=n_{i}$ at the point $\lambda^{\prime}$, take the limit, and conclude
that since $g_{n_{i}}(\lambda^{\prime})$ converges to $a$, $f_{n_{i}}(\lambda^{\prime} )$ converges to $1-a$. The preceding relation then shows  $a=1-a$, or $a=1/2$.
Since we have obtained convergence to 1/2 for $f_{n_{i}}(\lambda)$ along the arbitrary convergent subsequence $n_{i}$, it follows that $f_{n}(\lambda)$ itself converges to 1/2. The proof of Theorem~\ref{w12} is complete.

\section{Proof that Condition (TU) is satisfied}
Theorem~\ref{snb2} will follow from Theorem~\ref{w12} if it is shown that Condition (TU) is satisfied. We do that in this section, and begin with some elementary observations about the relationship between the digits in the expansion of $\omega$ and those in the expansion of $\nu=\omega^{2}$. 

By an {\it initial segment} of length $r$ of a dyadic expansion, we
refer to the
string of the first $r$ digits of the expansion. Let $\omega_{n}$ be the dyadic rational formed by the initial segment of length $n$ of $\omega$, that is,
$\omega_{n}=.x_{1}x_{2} \cdots x_{n}$.

For fixed $r>1$, consider the decomposition of $\Omega$ by the intervals
\[ I_{k+1}=[k\,2^{-r}, (k+1)\,2^{-r}) \hspace{.4in} 0 \leq k \leq
2^{r}-1. \]
Each of these intervals will be called an {\it r box}.
\begin{lemma} \label{inseg}
(a): If $\omega^{2} \in \,I_{k+1}$, then $.u_{1}u_{2} \cdots u_{r}$,
the initial segment of length $r$ of $\omega^{2}$, is equal to
$k\,2^{-r}$.\\
(b): Let $n$ digits $x_{1},x_{2}, \cdots ,x_{n}$ be specified. If
\[(.x_{1}x_{2} \cdots x_{n})^{2}  \mbox{  and  }
(.x_{1}x_{2} \cdots x_{n} + 2^{-n})^{2} \]
lie in the same $r$ box, say $I_{k+1}$, then, no matter how the
coordinates \\$x_{n+1},x_{n+2}, \cdots$ are subsequently chosen, the point
\begin{equation} \label{point}
\omega=.x_{1}x_{2} \cdots x_{n}x_{n+1} \cdots \end{equation}
is such that the initial segments of length $r$ of $\omega^{2}$ and of
each $\omega_{m}^{2}$ for $m \geq n$ are the same, with common value
$k\,2^{-r}$. \\
(c): Let $\omega^{2}$ be irrational. Given a positive integer $r$,
there exists a positive integer $N_{r}(\omega)>1$ such
that the initial segments of length $r$ of $\omega^{2}$ and of
each $\omega_{m}^{2}$ for $m \geq N_{r}$ are the same. Consequently,
each of the digits $u_{1}, u_{2}, \cdots, u_{r}$ is a function of the
digits $x_{m}, m \leq N_{r}$. Moreover, the set
$\{N_{r}(\omega)=n\}$ is defined in terms of
the coordinates $x_{1}, x_{2}, \cdots, x_{n-1}$ of $\omega$. \\
\end{lemma}
Proof: (a): The possible values of $.u_{1}u_{2} \cdots u_{r}$ are
\[.00\cdots00=0\cdot2^{-r},\, .00\cdots01=1\cdot2^{-r},\, .00\cdots10=2\cdot2^{-r},
\cdots \]
where a digit can change only if an amount at least equal to $2^{-r}$ is added
onto the current value. Therefore each
fixed value of $.u_{1}u_{2} \cdots u_{r}$ represents the left-hand
endpoint of a unique $r$ box. If $\omega^{2} \in\,I_{k+1}$, then
$\omega^{2}$ and $k\,2^{-r}$ differ by less than $2^{-r}$, so the initial
segment of $\omega^{2}$ is still $k\,2^{-r}$.\\
(b): The distance between the point
\[\omega_{1}=.x_{1}x_{2} \cdots x_{n} \mbox{  and  }
\omega_{2}=.x_{1}x_{2} \cdots x_{n}11\cdots \]
obtained by
choosing 1 for each $x_{m}, m >n$ is $2^{-n}$. Since the point $\omega$
of relation~\ref{point} satisfies $\omega_{1} \leq \omega \leq
\omega_{2}$, it follows from the assumption of part (b) that $\omega^{2}$
lies in $I_{k+1}$. The same argument holds for $\omega_{m}$ for $m \geq n$.
The conclusion now follows from part (a).\\
(c): Since $\omega^{2}$ is irrational, $\omega^{2}$ lies in the interior
of an $r$ box. As $m \rightarrow \infty$, $\omega_{m}^{2}$ tends to
$\omega^{2}$
and the $\omega_{m}^{2}$  are bounded away from the right-hand endpoint of
the $r$ box. It follows that eventually $\omega_{m}^{2}$  and
$(\omega_{m}+2^{-m})^{2}$ both lie in the same $r$ box. Define
$N_{r}(\omega)=n$ if $n$ is the smallest integer larger than 1 such that
$\omega_{n-1}^{2}$ and $(\omega_{n-1}+2^{-(n-1)})^{2}$ lie in the same
$r$ box. The assertion
about initial segments follows from part (b). From the definition of
$N_{r}$, the digits $u_{m}, m\leq r$ are determined by giving the first
$N_{r}-1$ $x$ coordinates, so the initial segment of length $r$ of the
$u$ sequence is a function of the initial segment of length $N_{r}-1$ of
the $x$ sequence. Moreover, to determine whether a point $\omega$ belongs to
$\{N_{r}(\omega)=n\}$, one need only know the first
$n-1$ coordinates of $\omega$. The last part of the assertion follows
readily from part (b). This completes the proof.

The preceding result showed that an initial segment of a $U$ sequence is determined by an initial segment of $X$ sequences. The reverse situation is also true: an initial segment of an $X$ sequence is determined by an initial segment of $U$ sequences. The argument is similar to that of Lemma~\ref{inseg}.
\begin{lemma} \label{inseg2}
Let $\nu$ be irrational, $\omega= \sqrt{\nu}$, and let $x_{1}, \cdots, x_{n}$ be the first $n$ coordinates of $\omega$. Then there exists an integer $m$ depending on $\nu$ and $n$ such that
if $\nu^{\prime}=u_{1}, \cdots u_{m},u^{\prime}_{m+1}, u^{\prime}_{m+2}, \cdots$ is any point whose initial $m$ segment agrees with that of $\nu$ but whose other coordinates may be arbitrary, then $\omega^{\prime}=\sqrt{\nu^{\prime}}$ has initial segment $x_{1}, \cdots, x_{n}$.
\end{lemma}
Proof:  The distance between $.u_{1}\cdots u_{j}$ and $.u_{1}\cdots u_{j}11\cdots$ is $2^{-j}$. Therefore the interval with endpoints $.u_{1}\cdots u_{j}$ and $.u_{1}\cdots u_{j}+2^{-j}$ contains $\nu^{\prime}$ and the endpoints converge to $\nu$. Decompose the unit interval into $n$ boxes (see Lemma~\ref{inseg}) and note that if $\omega_{1}$ and $\omega_{2}$ lie in the same $n$ box, the initial segments of length $n$ of each are the same. Since $\omega$ is irrational, it lies in the interior of an $n$ box. It follows that the square roots of the endpoints of the interval containing $\nu^{\prime}$ must eventually, for all sufficiently large $j$, be in the same $n$ box as $\omega$. Thus $\sqrt{\nu^{\prime}}$ must also be in this $n$ box. The proof is complete.

The following arguments will use some elementary ideas from the calculus of finite differences (see, e.g., \cite{sg:de}). We review some of the notation. Let $v(y_{1}, \cdots, y_{l})=v(\mbox{\boldmath $y$})$ be a function on the $l$-fold product space $S^{l}$ where the $y_{i} \in S$, a set of real numbers. 
Suppose that the variable $y_{i}$ is changed by the amount $\Delta y_{i}$ such that the $l$-tuple $\mbox{\boldmath $y^{(1)}$}=(y_{1}, \cdots, y_{l})$ is taken into $\mbox{\boldmath $y^{(2)}$}=(y_{1}+\Delta y_{1}, \cdots, y_{l}+\Delta y_{l})$ in the domain of definition of $v$. Put $v(\mbox{\boldmath $y^{(2)}$})-
v(\mbox{\boldmath $y^{(1)}$})= \Delta v$, and let 
\begin{eqnarray} \label{decomp1}
 \Delta v_{i}&=& v(y_{1}, \cdots, y_{i-1},y_{i}+\Delta y_{i}, y_{i+1}+\Delta y_{i+1, \cdots},
y_{l}+\Delta y_{l}) \hspace{.25in}\\
&-&v(y_{1}, \cdots, y_{i-1},y_{i}, y_{i+1}+\Delta y_{i+1, \cdots},
y_{l}+\Delta y_{l}). \nonumber \end{eqnarray} 
Then $\Delta v= \sum_{i} \Delta v_{i}$ is the total change in $v$ induced by changing all of the $y_{i}$, where this total change is written as a sum of step-by-step changes in the individual   
$y_{i}$. Formally, by dividing, we can write \begin{equation} \label{decomp2}
\Delta v= \sum_{i} (\Delta v_{i}/\Delta y_{i})\cdot \Delta y_{i}. \end{equation}
If some $\Delta y_{i}=0$, its coefficient in relation \ref{decomp2} has the form 0/0. Interpreting this coefficient as 0 makes the relation meaningful and true. Define the {\it partial difference of $v$ with respect to $y_{i}$, evaluated at the pair (\mbox{\boldmath $y^{(1)}$},\,\mbox{\boldmath $y^{(2)}$})} by 
\[\frac{\Delta v}{\Delta y_{i}}= \Delta v_{i}/ \Delta y_{i}.\]
Notice that the forward slash (/) in this relation expresses division and the horizontal slash on the left hand side is the partial difference operator. The sum of relation~\ref{decomp2} is called the {\it total difference of $v$} evaluated at the given pair and can also be written
\begin{equation} \label{decomp3}
\Delta v= \sum_{i} \frac{\Delta v}{\Delta y_{i}}\cdot \Delta y_{i}. \end{equation}
The partial and total differences are the discrete analogs of  the partial derivative of $v$ with respect to $y_{i}$ and the total differential,  respectively, in the theory of differentiable functions of several real variables. The partial difference of $v$ with respect to $y_{i}$ at a given pair is a measure of the contribution of $\Delta y_{i}$ to $\Delta v$ when all the other $y$ variables are held constant.   

We will say that $\omega$ and $\nu=\omega^{2}$ are points that {\it correspond} to one another. Since corresponding points uniquely determine each other, each $x_{i}$ is a function of the $u_{i}$ and the average $f_{n}(\omega)$ of relation~\ref{aver} can be written as a function $h_{n}(\nu)$ (see the Introduction). The function $f_{n}$ only depends on the first $n$ coordinates of $\omega$, and using a slight abuse of notation we understand by $f_{n}(x_{1}, \cdots, x_{n})$ the function of $n$ variables such that
\begin{equation} \label{xy2}
f_{n}(x_{1}, \cdots, x_{n})= f_{n}(\omega )=h_{n}(\nu)=h_{n}(u_{1}, u_{2}, \cdots). \end{equation}
Fix a point $\omega$ and for each $x_{i}$ let $\Delta x_{i}$ be a given increment ($\Delta x_{i}=0, 1,$ or $-1$). Let $\omega^{(1)}$ have coordinates $x_{i}+\Delta x_{i}$.  The changes $\Delta x_{i}$ correspond to changes $\Delta u_{i}$ in the coordinates of $\nu$, the point corresponding to $\omega$, such that  $\nu$ goes into the point $\nu^{(1)}$ with coordinates $u_{i}+\Delta u_{i}$ corresponding to $\omega^{(1)}$. Assume that all $X$ and $U$ sequences discussed here and below represent irrational numbers. Now consider the change 
\[f_{n}(\omega^{(1)})-f_{n}(\omega)=\Delta f_{n}=\Delta h_{n}=h_{n}(\nu^{(1)})-h_{n}(\nu),\]
where the right hand side can be written \begin{equation} \label{aitch}  
\Delta h_{n}= h_{n}(u_{1}+\Delta u_{1}, u_{2}+\Delta u_{2}, \cdots)-h_{n}(u_{1}, u_{2} \cdots).\end{equation}

Recall that the capital letter notation $X_{i}$ and $U_{i}$ denotes the $i$th coordinate variable of $\omega$ and $\nu$, respectively. This notation will be convenient when small letters may be reserved to denote particular values.
\begin{lemma} \label{fini} At the pair $(\nu, \nu^{(1)})$, $\Delta h_{n}$ can be represented as a total difference
\begin{eqnarray} \label{delh}
\Delta h_{n}&=& h_{n}(u_{1}+\Delta u_{1}, u_{2}+\Delta u_{2}, \cdots)-h_{n}(u_{1}, u_{2}, \cdots)\nonumber \\& & = \sum_{i \geq 1} \frac{\Delta h_{n}}{\Delta U_{i}}\, \Delta U_{i}=          \sum_{i \geq 1} (\Delta h_{n,i}/\Delta U_{i})\, \Delta U_{i}, \end{eqnarray}  
where  $\Delta U_{i}= \Delta u_{i}$ and
\begin{eqnarray} \label{pd}  \Delta h_{n,i}=\Delta h_{n,i}(\nu, \nu^{(1)})=&&\\
h_{n}(u_{1}, \cdots, u_{i-1}, u_{i}+\Delta u_{i}, u_{i+1}+\Delta u_{i+1}, \dots)-&& \nonumber\\
h_{n}(u_{1}, \cdots, u_{i-1}, u_{i}, u_{i+1}+\Delta u_{i+1}, \dots)&&. \nonumber \end{eqnarray}
The formally infinite sum of relation~\ref{delh} reduces to a finite sum when evaluated at the  pair $(\nu, \nu^{(1)})$, that is, given the pair, there exists an integer $m$ such that the partial differences  $\Delta h_{n,i}/\Delta U_{i}=0$ for all $i>m$. The number of non-vanishing terms in the sum depends on the pair chosen and on $n$.
\end{lemma}
Proof: The function $h_{n}=f_{n}$ only depends on the initial segment of length $n$ of $\omega$. Given $\nu$, Lemma~\ref{inseg2} proves the existence of an integer $m$ such that for all $i>m$, the points with coordinates 
\begin{eqnarray}
u_{1}, \cdots, u_{i-1}, u_{i}+\Delta u_{i}, u_{i+1}+\Delta u_{i+1}, \dots &\mbox{and}&\\
 u_{1}, \cdots, u_{i-1}, u_{i}, u_{i+1}+\Delta u_{i+1}, \dots&& \nonumber \end{eqnarray}
correspond to $X$ sequences having the same initial segment of length $n$ as $\omega$. Consequently, for $i>m$ the difference terms in relation~\ref{pd} are equal and the partial differences evaluated at the given pair vanish. The argument is thus reduced to the observations leading to relations~\ref{decomp1} and \ref{decomp2}, and the proof is finished. 

Since $X_{j}$ is a function of the $U$ variables for all $j$, an argument similar to the above shows that there is a representation analogous to relation~\ref{delh} of the form
\begin{equation} \label{delh2}
\Delta X_{j}=\sum_{i \geq 1}\frac{\Delta X_{j}}{\Delta U_{i}}\, \Delta U_{i}=\sum_{i \geq 1}
(\Delta X_{j,i}/\Delta U_{i})\, \Delta U_{i}, \end{equation}  
where $\Delta X_{j,i}$ is derived from $\Delta X_{j}$ in the same way as $\Delta h_{n,i}$ is derived from $\Delta h_{n}$ (see relation~\ref{pd}).
At a given pair this representation also reduces to a finite sum by Lemma~\ref{inseg2}.

In similar fashion, representations in terms of the $X$ variables may be written. At the pair $(\omega, \omega^{(1)})$ we obtain the following relations analogous to relations~\ref{delh} and \ref{delh2}:    
\begin{equation} \label{delf}
\Delta h_{n}=\Delta f_{n}=\sum_{j \geq 1}
\frac{\Delta f_{n}}{\Delta X_{j}}\, \Delta X_{j}=\sum_{j \leq n}
\frac{1}{n}\, \Delta X_{j},  \end{equation}  
and
\begin{equation} \label{delf2}
\Delta U_{i}=\sum_{j \geq 1}\frac{\Delta U_{i}}{\Delta X_{j}}\, \Delta X_{j}=\sum_{j \geq 1}
(\Delta U_{i,j}/\Delta X_{j})\, \Delta X_{j}, \end{equation}  
where $\Delta U_{i,j}$ is defined similarly to $\Delta h_{n,i}$ and $\Delta X_{j,i}$.

Now fix the positive integer $k$ and refer to relation~\ref{delh}. It will be seen
that at the pair $(\nu, \nu^{(1)})$ and for any subsequence $n_{l}$, $\limsup_{l} \Delta h_{n_{l}}$ does not depend on the values of $\Delta U_{i},\, i \leq k$. This is what must be shown to prove the validity of Condition (TU), that is, that the influence of any initial segment of the $U$ variables on $h_{n}$ dies out in the limit. 
To this end we first observe that either the $X$ or the $U$ variables may be taken as independent variables, with the other set dependent on them. We take the $X$ variables as independent. We now rewrite  relation~\ref{delh} in terms of the $X$ variables, using the dependence relation of $U$ on $X$. The following result is an analog of the chain rule for differentiable functions of several real variables.
\begin{lemma}  \label{indyx} Let the $U$ variables be functions of the independent $X$ variables. Then for each $j$ the following relations are valid evaluated at any pair: 
 \begin{equation} \label{chrule}
 \Delta h_{n}= \sum_{j} \left(\sum_{i}  \frac{\Delta h_{n}}{\Delta U_{i}}\frac{\Delta U_{i}}{\Delta X_{j}}\right) \Delta X_{j}, \end{equation}
and
\begin{equation} \label{chrule2}
\sum_{i} \frac{\Delta h_{n}}{\Delta U_{i}}\frac{\Delta U_{i}}{\Delta X_{j}}=\frac{\Delta h_{n}}{\Delta X_{j}}=\frac{1}{n}\,\, \mbox{ for } j \leq n \mbox{ and $=0$ for } j>n.  \end{equation}
\end{lemma}
Proof: Start with relation~\ref{delh}, substitute relation~\ref{delf2} and interchange the order of addition (possible because of the finiteness of the sums) to get: 
\begin{eqnarray}\Delta h_{n}=\sum_{i} \frac{\Delta h_{n}}{\Delta U_{i}} \Delta U_{i} &=&\sum_{i} \frac{\Delta h_{n}}{\Delta U_{i}} \left(\sum_{j}\frac{\Delta U_{i}}{\Delta X_{j}}\Delta X_{j}\right)\\ &=&\sum_{j}\left(\sum_{i}\frac{\Delta h_{n}}{\Delta U_{i}}\frac{\Delta U_{i}}{\Delta X_{j}}\right) \Delta X_{j}, \nonumber \end{eqnarray}
and this is relation~\ref{chrule}. To prove relation~\ref{chrule2}, take the partial difference on both sides of relation~\ref{chrule} with respect to $X_{j_{0}}$ for a fixed index $j_{0}$. Note that the partial difference of a sum evaluated at a given pair is additive, so that 
\begin{equation} \label{j0}
\frac{\Delta h_{n}}{\Delta X_{j_{0}}}= \sum_{j} \left(\sum_{i}  \frac{\Delta h_{n}}{\Delta U_{i}}\frac{\Delta U_{i}}{\Delta X_{j}}\right)\frac{ \Delta X_{j}}{\Delta X_{j_{0}}}, \end{equation}
where the left hand side of relation~\ref{j0} is equal to $1/n$ or 0 depending on whether $j_{0} \leq n$ or $>n$. The variables $X_{j}$ are independent. This means that a change in $X_{j_{0}}$ does not cause a change in $X_{j}, j \neq j_{0}$, that is
 \[\frac{\Delta X_{j}}{\Delta X_{j_{0}}}=0, j \neq j_{0} \mbox{ and } =1 \mbox{ if } j=j_{0}. \] 
Relation~\ref{j0} thus reduces to relation~\ref{chrule2}.
\begin{lemma} \label{ndk}
Let $k$ be a fixed index and let the pair $(\nu, \nu^{(1)})$ be given as in Lemma~\ref{fini}. Then at any pair $(\nu, \nu^{(*)})$ for which $\nu^{(*)}$ has the same tail as $\nu^{(1)}$ starting from index $(k+1)$, we have 
\begin{eqnarray} \label{delhk}
\limsup_{n} \Delta h_{n} =\limsup_{n} \sum_{i>k}
\frac{\Delta h_{n}}{\Delta U_{i}} \Delta U_{i}. \end{eqnarray}
The right hand side of relation~\ref{delhk} is the same, term by term, as the right hand side of relation~\ref{delh} if $\Delta U_{i},\, i \leq k$ are set equal to 0 there.   
Thus $\limsup_{n} \Delta h_{n}$ is constant over all such pairs $(\nu, \nu^{(*)})$ and does not depend on the values of $\Delta U_{i}, \, i \leq k$. Relation~\ref{delhk} remains true if the limit superior is taken over any subsequence $n_{l}$ rather than the entire sequence. In addition, in relation~\ref{delh} we have
 \[\lim_{n}\Delta h_{n,k}= \lim_{n}\frac{\Delta h_{n}}{\Delta U_{k}}=0.\]  
 \end{lemma}
Proof: Let $k$ be given, and consider the set $S_{1}$ of the $2^{k}$ points in $U$ space for which the tail starting from the coordinate $(k+1)$ is the same as that of $\nu^{(1)}$ but the initial segment of length $k$ runs through all possible values. Let $S$ be the set of the $2^{k}$ pairs $(\nu, \nu^{(*)})$ where $\nu^{(*)} \in S_{1}$. According to Lemma~\ref{inseg} there exists $N$ such that $U_{i},\,\, i \leq k$ is a function of the corresponding $X_{j},\,\, j \leq N$ for all points in $S_{1}$. Relation~\ref{chrule2} implies that at any pair whatsoever
\begin{equation} \label{init3}
 \lim_{n} \sum_{j \leq N} \left(\sum_{i}  \frac{\Delta h_{n}}{\Delta U_{i}}\frac{\Delta U_{i}}{\Delta X_{j}}\right) \Delta X_{j} =0. \end{equation}
 Relation~\ref{chrule} then gives   \begin{eqnarray} \label{2terms}
 \limsup_{n} \Delta h_{n}= \limsup_{n} \Delta f_{n}=\limsup_{n} \sum_{j>N} \left(\sum_{i}  \frac{\Delta h_{n}}{\Delta U_{i}}\frac{\Delta U_{i}}{\Delta X_{j}}\right)  \Delta X_{j}= & &  \\ 
 \limsup_{n} \left[\sum_{j>N} \left(\sum_{i\leq k}  \frac{\Delta h_{n}}{\Delta U_{i}}\frac{\Delta U_{i}}{\Delta X_{j}}\right) \Delta X_{j}
+\sum_{j>N} \left(\sum_{i>k}  \frac{\Delta h_{n}}{\Delta U_{i}}\frac{\Delta U_{i}}{\Delta X_{j}}\right) \Delta X_{j}\right].& & \nonumber \end{eqnarray}

We claim that at pairs of $S$ the first of the two double sums within the bracket on the right hand side of relation~\ref{2terms} vanishes. The reason is that $U_{i}, i\leq k$ only depend on $X_{j}, j \leq N$ so that 
\[\frac{\Delta U_{i}}{\Delta X_{j}}=0 \mbox{ for } i \leq k \mbox{ and } j>N. \] 
It follows that we have  \begin{equation} \label{final}
\limsup_{n} \Delta f_{n}=\limsup_{n}\sum_{j>N} \left(\sum_{i>k}  \frac{\Delta h_{n}}{\Delta U_{i}}\frac{\Delta U_{i}}{\Delta X_{j}}\right) \Delta X_{j}. \end{equation}
It should be noted that the partial differences of $h_{n}$ with respect to $U_{i}$
for $i>k$ in relation~\ref{final} at a pair of $S$ are the same as the corresponding terms in relation~\ref{delh} at the pair $(\nu, \nu^{(1)})$; this is immediate by comparison using relation~\ref{pd}. Also observe
that relations~\ref{2terms} and \ref{final} also hold if the limit superior had been taken over any subsequence $n_{l}$ instead of over the entire sequence.

Let us define
\[h_{n}^{(k)}(U_{k+1}, U_{k+2}, \cdots)=h_{n}(u_{1},\cdots, u_{k}, U_{k+1}, U_{k+2}, \cdots)\]
where we recall that $u_{1},\cdots, u_{k}$ is the initial segment of length $k$ of $\nu$ and $U_{k+1}, U_{k+2}, \cdots$ are variables. Then at any pair in $S$ the value of $\Delta h_{n}^{(k)}$ is constant, namely \begin{equation} \label{okey}
\Delta h_{n}^{(k)}=\sum_{i>k }\frac{\Delta h_{n}}{\Delta U_{i}} \Delta U_{i}, \end{equation}
where again we observe that the terms are the same as the corresponding terms in relation~\ref{delh}. Now for any pair in $S$ we may also write a representation of $\Delta h_{n}^{(k)}$ in terms of the $X$ variables:
\begin{equation} \label{hnkxj}
 \Delta h_{n}^{(k)}= \sum_{j \geq 1} \frac{\Delta h_{n}^{(k)}}{\Delta X_{j}} \Delta X_{j}. \end{equation}
For $i \leq k$, the partial difference of $\Delta h_{n}^{(k)}$ with respect to $\Delta U_{i}$ is equal to 0, and so  for any $j$ \begin{equation} \label{lst1}
 \frac{\Delta h_{n}^{(k)}}{\Delta X_{j}}= \sum_{i \geq 1}  \frac{\Delta h_{n}^{(k)}}{\Delta U_{i}}\frac{\Delta U_{i}}{\Delta X_{j}}=\sum_{i>k}  \frac{\Delta h_{n}^{(k)}}{\Delta U_{i}}\frac{\Delta U_{i}}{\Delta X_{j}}=\sum_{i>k }  \frac{\Delta h_{n}}{\Delta U_{i}}\frac{\Delta U_{i}}{\Delta X_{j}}. \end{equation}
Put relation~\ref{lst1} into relation~\ref{hnkxj} to get \begin{equation} \label{same}
\limsup_{n} \Delta h_{n}^{(k)}=\limsup_{n}\sum_{j \geq 1} \left(\sum_{i>k}  \frac{\Delta h_{n}}{\Delta U_{i}}\frac{\Delta U_{i}}{\Delta X_{j}}\right) \Delta X_{j}. \end{equation}
The argument leading to  relation~\ref{same} goes through without change for any subsequence $n_{l}$ used instead of the entire sequence. The left hand side of relation~\ref{same} does not depend on any fixed coordinate $X_{j_{0}}$, even if the limsup is taken over any subsequence.  
If the coefficients of $\Delta X_{j_{0}}$ in relation~\ref{same}  \begin{equation} \label{nodep}
 \sum_{i>k}  \frac{\Delta h_{n}}{\Delta U_{i}}\frac{\Delta U_{i}}{\Delta X_{j_{0}}} \end{equation}
do not converge to 0, then there is a subsequence converging to a non-zero value. This would mean that the limsup of the right hand side of relation~\ref{same} along this subsequence would depend on the value of $\Delta X_{j_{0}}$. But this contradicts what is known for the left hand side of this relation. So the terms in relation~\ref{nodep} converge to 0 and this holds for arbitrary indices $j_{0}$, permitting us to rewrite relation~\ref{same} as 
 \begin{equation} \label{same2}
\limsup_{n} \Delta h_{n}^{(k)}=\limsup_{n}\sum_{j>N} \left(\sum_{i>k}  \frac{\Delta h_{n}}{\Delta U_{i}}\frac{\Delta U_{i}}{\Delta X_{j}}\right) \Delta X_{j}. \end{equation}
The right hand sides of relations~\ref{same2} and \ref{final} are the same. These relations together with relation~\ref{okey} prove that at any pair of $S$
\begin{equation} \label{tail2}
\limsup_{n} \Delta f_{n}=\limsup_{n} \Delta h_{n}=\limsup_{n}\Delta h_{n}^{(k)}=\limsup_{n} \sum_{i>k}  \frac{\Delta h_{n}}{\Delta U_{i}}\,\,\Delta U_{i}. \end{equation}
The argument leading to  relation~\ref{tail2} goes through without change for any subsequence $n_{l}$ used instead of the entire sequence. This relation thus expresses $\limsup_{l} \Delta h_{n_{l}}$  for an arbitrary subsequence at any pair of $S$ in terms of a function of the $U_{i}$ variables for $i>k$ which is constant at pairs in $S$. This function is the same as that in relation~\ref{delh} if $\Delta U_{i}=0,\,\, i \leq k$. It will follow that $\lim_{n}\Delta h_{n,k}=0$ in relation~\ref{delh}; the argument is similar to one given above. 
Without loss of generality take $k=1$ and consider the representation of $\Delta h_{n}$ at $(\nu, \nu^{(1)})$ given by relation~\ref{delh}. If there is a subsequence  $\Delta h_{n_{l},1}$, say, converging, if possible, to a non-zero value, then relation~\ref{delh} proves $\limsup_{l} \Delta h_{n_{l}}$ dependent on $\Delta U_{1}$, but this contradicts relation~\ref{tail2} for subsequences. Therefore $\Delta h_{n,1}$ converges to 0, and the proof of the lemma is complete. 
   
We are ready to summarize the foregoing results into a formal statement that Condition (TU) is satisfied.
\begin{lemma} \label{tu} (Condition (TU) is satisfied)
Let $f_{n}(\omega)$ be the average of relation~\ref{aver}. Let  $n_{i}$  be any fixed subsequence. 
Then the function $f=\limsup_{i \rightarrow \infty}f_{n_{i}}$ is a tail function with respect to the variables $U_{1}, U_{2}, \cdots$, that is, for any given positive integer $k$, the function $f$ can be written as a function of $U_{k+1}, U_{k+2}, \cdots$. Consequently, if 
$\omega$ and $\omega^{(1)}$ have corresponding points $\nu$ and $\nu^{(1)}$ with the same tail, then $f(\omega)=f(\omega^{(1)})$.\end{lemma}
Proof: Suppose that $\nu$ and $\nu^{(1)}$ have the same tail starting from index $k+1$. Apply Lemma~\ref{ndk} to conclude that $\limsup_{n} \Delta h_{n}=\limsup_{n} \Delta f_{n}$ only depends on the values of $\Delta U_{i}$ for $i>k$ and has the representation given by relation~\ref{delhk}.  
This proves the assertion and concludes the proof of Theorem~\ref{snb2}. 

We end our discussion by mentioning the problem of extending Theorem~\ref{snb2} from simple normality to normality. It appears that an approach similar to the one given here will work. We hope to have completed results soon. 


\newpage

\begin{thebibliography}{99}
\bibitem{bc:nn} David H. Bailey and Richard E. Crandall, {\it Random Generators and Normal Numbers}, to appear in Experimental Mathematics (also see first author's web site).
\bibitem{sg:de} Samuel Goldberg, {\it Introduction to Difference Equations}, John Wiley, 1958.
\bibitem{gh:pm} G. H. Hardy, {\it A Course of Pure Mathematics}, Cambridge University Press, 1952.
\bibitem{ni:ir} Ivan Niven, {\it Irrational Numbers}, Mathematical
Association of America (Carus Mathematical Monograph), distributed by
John Wiley, 1956.
\end{thebibliography}
\end{document}